\documentclass[reqno]{amsart}
\usepackage{amsmath, amsfonts, epsfig, amssymb, graphics, hyperref, mathrsfs}
\usepackage{rotating} 
\usepackage[usenames]{color}

{\theoremstyle{plain}
\newtheorem{theorem}{Theorem}[section]

\newtheorem{proposition}[theorem]{Proposition}
\newtheorem{lemma}[theorem]{Lemma}
\newtheorem{corollary}[theorem]{Corollary}
}
{\theoremstyle{definition}
\newtheorem{definition}[theorem]{Definition}
\newtheorem{remark}[theorem]{Remark}
\newtheorem*{note}{Note}
\newtheorem{example}[theorem]{Example}
}

\numberwithin{equation}{section}

\newlength{\raisesize} 
\makeatletter 
\ifthenelse{\@ptsize = 0}{\setlength{\raisesize}{5pt}}{}
\ifthenelse{\@ptsize = 1}{\setlength{\raisesize}{5.5pt}}{}
\ifthenelse{\@ptsize = 2}{\setlength{\raisesize}{6pt}}{}
\makeatother 

\newcommand{\des}[1]{D(#1)}
\newcommand{\maj}[1]{\mathrm{maj}(#1)}
\newcommand{\gmaj}[2]{\mathrm{maj}_{#1}(#2)}
\newcommand{\gmaja}[1]{\gmaj{\al}{#1}}
\newcommand{\sn}{S_n}
\newcommand{\al}{\mathbf{q}}
\newcommand{\ali}{q}
\newcommand{\K}{K}
\newcommand{\e}[2]{\kappa_{#1}(#2)}
\newcommand{\shuffle}{\raisebox{\raisesize}{\mbox{\:\begin{turn}{270} $\exists$ \end{turn}\:}}}
\newcommand{\kx}{K \langle X \rangle}
\newcommand{\scalar}[2]{\langle #1, #2 \rangle}

\newcommand{\p}{p}
\newcommand{\puv}{P_{u,v}}
\newcommand{\linext}[1]{L(#1)}
\newcommand{\apw}{\mathcal{A}(P,\omega)}
\newcommand{\fpwx}{F(P,\omega; \x)}
\newcommand{\id}{\text{id}}
\newcommand{\twisted}{\ltimes}
\newcommand{\ksn}{\K \sn}
\newcommand{\genk}{\K(\al)}
\newcommand{\genksn}{\K(\al) \sn}
\newcommand{\ru}{\zeta}
\newcommand{\partner}[2]{\eta_{#1}(#2)}
\newcommand{\kly}{\kappa_n}
\newcommand{\klypartner}{\eta_n}
\newcommand{\cycle}{\gamma}
\newcommand{\lie}{\e{n}{\al}}
\newcommand{\prt}{\partner{n}{\al}}
\newcommand{\N}[1]{N_{\al}(#1)}
\newcommand{\D}[1]{D_{\al}(#1)}
\newcommand{\apply}[2]{#1[#2]}
\newcommand{\liealg}{\mathcal{L}_n(\al)}
\newcommand{\liealgn}{\mathcal{L}_n}
\newcommand{\idem}{\pi}
\newcommand{\freelie}{\mathcal{L}_K(X)}
\newcommand{\xj}{x}
\newcommand{\x}{\mathbf{x}}
\newcommand{\infprod}{\Theta(\x)}
\newcommand{\infproda}{\Theta(\al)}
\newcommand{\kxn}{\K[[\xj_1, \ldots, \xj_n]]}
\newcommand{\directsum}{\oplus_{n\geq 0}\kxn \sn}
\newcommand{\mare}{\star}
\newcommand{\st}{\mathrm{st}}
\newcommand{\emptyperm}{\varepsilon}
\newcommand{\fancyprod}[1]{\overleftarrow{\prod_{#1}}}

\begin{document}
\title[A Multi-Parameter Klyachko Idempotent]{
$P$-partitions and a multi-parameter Klyachko idempotent}

\author{Peter McNamara}
\address{Instituto Superior T\'ecnico\\
Departamento de Matem\'atica\\
Avenida Rovisco Pais\\
1049-001 Lisboa\\
Portugal.}
\email{mcnamara@math.ist.utl.pt}

\author{Christophe Reutenauer}
\address{Laboratoire de Combinatoire et d'Informatique Math\'ematique\\
Universit\'e du Qu\'ebec \`a Montr\'eal\\
Case Postale 8888, succursale Centre-ville\\
Montr\'eal (Qu\'ebec) H3C 3P8\\
Canada.}
\email{christo@lacim.uqam.ca}

\dedicatory{Dedicated to Richard Stanley on the occasion of his 60th birthday}

\begin{abstract}
Because they play a role in our understanding of
the symmetric group algebra, Lie idempotents
have received considerable attention.
The Klyachko idempotent has attracted interest from combinatorialists,
partly because its definition involves the major index of permutations.

For the symmetric group $\sn$,
we look at the symmetric group algebra with coefficients from the
field of rational functions in $n$ variables $\ali_1, \ldots, \ali_n$.
In this setting, we can define an $n$-parameter generalization
of the Klyachko idempotent, and we show it is a Lie idempotent
in the appropriate sense.  Somewhat surprisingly, our proof that it is a Lie element
emerges from Stanley's theory of $P$-partitions.
\end{abstract}

\maketitle
  

\section{Introduction}

The motivation for our work is centered around the search for Lie idempotents in the
symmetric group algebra.  In fact, our goal is to give a generalization of the
well-known Klyachko idempotent, and to show that important and interesting properties of
the Klyachko idempotent carry over to the extended setting.  It turns out that the proof
that our generalized Klyachko idempotent is a Lie element
gives a nice application and illustration of Richard Stanley's theory of $P$-partitions.
We should point out that $P$-partitions were previously used in \cite{BBG90} to 
show that the traditional Klyachko idempotent is a Lie element.  

To define Lie idempotents,
however, we will first need the concepts of free Lie algebras and the symmetric group
algebra.
Let $\K$ be a field of characteristic $0$.  If $X$ is an alphabet, we will write $\kx$ to
denote the free associative algebra consisting of all linear combinations of words
on $X$ with coefficients in $\K$.
The product of two words on $X$ is defined to be their concatenation, and extending
this product by linearity gives a product on $\kx$.  We can then define the 
\emph{Lie bracket} $[p, q]$ of two elements $p$ and $q$ of $\kx$ by $[p,q] = pq - qp$. 
We let $\freelie$ denote the smallest vector subspace of $\kx$ containing $X$ and closed
under the Lie bracket.  It is a classical result that $\freelie$ is the free Lie algebra
on $X$. We refer the reader to \cite{Gar90,Reu93} for further details on free 
Lie algebras from a combinatorial viewpoint. 

If $X = \{1,2,\ldots, n\}$, then elements of the symmetric group $\sn$ can be considered as
words on $X$.   We write $\ksn$ to denote 
the \emph{symmetric group algebra}, which consists of linear combinations of elements
of the symmetric group $\sn$, with coefficients in $\K$.  
When $X = \{1,2,\ldots, n\}$, certain elements of $\freelie$, such as 
$[[\cdots[1,2],3],\ldots,n]$, can be naturally
considered to be elements of $\ksn$.  This is because all the words in their expansions as  
elements of $\kx$ are permutations of $1,2,\ldots,n$.  
Elements in this intersection of $\ksn$ and $\freelie$ are called \emph{Lie elements}.
We will denote the set of Lie elements by $\liealgn$.  

We should clarify our suggestion that $\ksn$ is an algebra.    
The product
of two permutations $\sigma, \tau \in \sn$ is the usual composition $\sigma \tau$ from 
right to left, and extending this product by linearity gives the product in $\ksn$.  
It is well-known, and is not difficult to check, that $\liealgn$ is then a left ideal of
$\ksn$.

\begin{definition}\label{def:lieidempotent}
A Lie idempotent is an element $\idem$ of $\ksn$ that is idempotent and that satisfies
\[
\ksn \idem = \liealgn.
\]
\end{definition}
In particular, $\idem$ must be a Lie element. 
Lie idempotents are quite remarkable 
because, in particular, they give an alternative, and direct, construction of $\liealgn$.
It is natural, therefore, that there should be widespread interest in the search for Lie idempotents
and, from a combinatorial perspective, \cite{BBG90,DKLT94,GaRe89,Gar90,KLT97,PaRe02}
all offer progress in this search.

One of the most famous Lie idempotents is the Klyachko idempotent of \cite{Kly74}.  (We refer 
the reader to the end of this introduction for the definition of the major index, $\maj{\sigma}$,
of $\sigma \in \sn$.)
Let $\ru$ be a primitive $n$th root of unity in $\K$, meaning that $\ru^n = 1$ and $\ru^m \neq 1$ for 
$1 \leq m < n$.  Then the Klyachko idempotent $\kly$ is defined by
\[
\kly = \frac{1}{n} \sum_{\sigma \in \sn} \ru^{\maj{\sigma}} \sigma . 
\]
The appearance of the major index in this definition naturally makes the Klyachko 
idempotent appealing to combinatorialists and, for example, \cite{BBG90,DKLT94,KLT97,PaRe02}
study the Klyachko idempotent and its generalizations.

Our goal is to introduce a new, broad generalization of the Klyachko idempotent and
to show that its Lie idempotency property is preserved in this much wider setting.
As we will show in Example \ref{exa:klyachko}, we will indeed be able to recover 
the usual symmetric group algebra $\ksn$ and the Klyachko idempotent by specialization.  
Rather than working with a primitive $n$th root of unity, we will let $\al = (\ali_1, \ali_2, \ldots, \ali_n)$
be a sequence of variables in a field $\K$ with the only restriction
being that $\ali_1 \ali_2 \cdots \ali_n = 1$.
Since the $\ali_i$'s are formal variables, we can assume, in particular, that
$\ali_{i_1} \ali_{i_2} \cdots \ali_{i_r} \neq 1$ for any 
proper subset $\{i_1, i_2, \ldots, i_r\}$ of $\{1, 2, \ldots, n\}$.
Throughout, unless otherwise stated, $\al$ will denote such a sequence.
We let $\genk$ denote the field of rational functions in $\al$ over the field $\K$, 
and our primary focus will 
be $\genksn$, the symmetric group algebra with coefficients in $\genk$.  
Before proceeding, however, we must pay attention to a twist in our story.
It turns out that the most useful product for $\genksn$ is not the natural 
analogue of the product for $\ksn$.  
More precisely, if $f(\al), g(\al) \in \genk$ and $\sigma, \tau \in \sn$, 
one might assume that the product of
$f(\al)\sigma$ and $g(\al)\tau$ should be defined to be simply $(f(\al)g(\al)) \sigma\tau$.  
However, this product
does not seem to allow the concepts of interest from $\ksn$ to extend to $\genksn$
and, in particular, our generalized Klyachko element is not idempotent with 
respect to this product for $n \geq 4$.
Instead, we observe that there is a natural left action of $\sn$ on $\genk$: if 
$f(\al) = f(\ali_1, \ali_2, \ldots, \ali_n) \in \genk$, then we define  
\[
\apply{\sigma}{f(\al)} = f(\ali_{\sigma(1)}, \ali_{\sigma(2)}, \ldots, \ali_{\sigma(n)}).
\]
We then define the \emph{twisted product} of $f(\al)\sigma$ and $g(\al)\tau$,
denoted $f(\al)\sigma \twisted g(\al)\tau$, by
\[
f(\al)\sigma \twisted g(\al)\tau = (f(\al) \apply{\sigma}{g(\al)}) \sigma\tau.
\]
As a simple example, if n=3,
\[
231 \twisted \frac{\ali_1 \ali_3}{(1-\ali_1)(1-\ali_1 \ali_3)} 132 = 
\frac{\ali_2 \ali_1}{(1-\ali_2)(1-\ali_2 \ali_1)} 213.
\]
This twisted product appears in some standard texts on the representation 
theory of groups and algebras, 
such as \cite[\S 28]{CuRe81}.
As in \cite{CuRe81}, 
we leave it as a quick exercise to check that the twisted product is associative.  
Our results will serve as evidence in favor of the assertion that the twisted product is
the ``correct'' product for $\genksn$.  

We are now in a position to define our extended version
of the Klyachko idempotent.

\begin{definition} 
Given a permutation $\sigma$ in $\sn$,
define the $\al$-major index $\gmaja{\sigma}$ of $\sigma$ by
\[
\gmaja{\sigma} = \frac{\prod_{j \in \des{\sigma}} \ali_{\sigma(1)} \ali_{\sigma(2)} \ldots \ali_{\sigma(j)}}
{\prod_{i=1}^{n-1} (1 - \ali_{\sigma(1)} \ali_{\sigma(2)} \ldots \ali_{\sigma(i)})}.
\]
\end{definition}

We will justify the terminology ``$\al$-major index'' in Example \ref{exa:klyachko}.

\begin{remark}
The numerator terms $\N{\sigma} = \prod_{j \in \des{\sigma}} \ali_{\sigma(1)} \ldots \ali_{\sigma(j)}$ 
have a certain fame due to their appearance in \cite{Gar80}.  There, 
Garsia shows that every polynomial $G(\al)$ in $\ali_1, \ldots, \ali_n$ has a unique expression
of the form
\[
G(\al) = \sum_{\sigma \in \sn} g_\sigma(\al) \N{\sigma},
\]
where each $g_\sigma(\al)$ is a polynomial that is \emph{symmetric}
in $\ali_1, \ldots, \ali_n$.  Furthermore, if $G(\al)$ has integer coefficients, then so
do all the polynomials $g_\sigma(\al)$.  These results were originally
conjectured by Ira Gessel.
\end{remark}

The following definition introduces our main object of study.

\begin{definition} 
Denote by $\lie$ the element of $\genksn$ given by
\[
\lie = \sum_{\sigma \in \sn} \gmaja{\sigma}\, \sigma .
\]
\end{definition}

\begin{example}
If $n=3$ we get
\begin{eqnarray*}
\e{3}{\al} & = & 
\frac{1}{(1-\ali_1)(1-\ali_1 \ali_2)} 123 +
\frac{\ali_1 \ali_3}{(1-\ali_1)(1-\ali_1 \ali_3)} 132 \\ 
& & + \frac{\ali_2}{(1-\ali_2)(1-\ali_1 \ali_2)} 213 
+ \frac{\ali_2 \ali_3}{(1-\ali_2)(1-\ali_2 \ali_3)} 231 \\
& & + \frac{\ali_3}{(1-\ali_3)(1-\ali_1 \ali_3)} 312 +
\frac{\ali_2 \ali_3^2}{(1-\ali_3)(1-\ali_2 \ali_3)} 321 .
\end{eqnarray*}
\end{example}

\begin{example}\label{exa:klyachko}
With $\ru$ a primitive $n$th root of unity, we can see that $\lie$ maps to the Klyachko
idempotent $\kly$ under
the specialization $\ali_i \to \ru$ for $i=1, \ldots, n$.
Indeed, the $\al$-major index of any $\sigma \in \sn$ then specializes to 
\[
\frac{\ru^{\maj{\sigma}}}{\prod_{i=1}^{n-1} (1 - \ru^i)} 
= \frac{\ru^{\maj{\sigma}}}{n}.
\]
This equality follows from the identity
\[
x^n-1 = (x-1)(x-\ru)(x-\ru^2)\cdots (x-\ru^{n-1}),
\]
which implies that $n = (1-\ru)(1-\ru^2)\cdots(1-\ru^{n-1})$.  
Also notice that, since $q_1 = \cdots = q_n$, the twisted product of $\genksn$ 
in this setting
is identical to the usual product of $\ksn$. 

Actually, if we take the ring $\K[q]$, localized at $1-q, \ldots, 1-q^{n-1}$, 
and quotiented by the ideal generated by $1-q^n$, it will follow from our results that the element
\[
\sum_{\sigma \in \sn} \frac{q^{\maj{\sigma}}}{(1-q)\cdots(1-q^{n-1})}\sigma
\]
is a Lie idempotent.  This may shed some light on some results in \cite{BBG90}.  
We see that this new element specializes to the Klyachko idempotent when we map
$q$ to a \emph{primitive} root of unity.  
\end{example}

We can now state our main results.  

\begin{theorem}\label{thm:lie} 
$\lie$ is a Lie element. 
\end{theorem}

\begin{theorem}\label{thm:idem}
$\lie \twisted \lie = \lie$, i.e., $\lie$ is idempotent as an element of $\genksn$.
\end{theorem}

We will let $\liealg$ denote the analogue of $\liealgn$ when coefficients come from $\genk$.  As
a formula, $\liealg = \genksn \cap \mathcal{L}_{\genk}(\{1,2,\ldots,n\})$.
We observe that $\liealg$ is a left ideal of $\genksn$:
\begin{equation}\label{equ:liealgleftideal}
\genksn \twisted \liealg = \liealg.
\end{equation}

\begin{theorem}\label{thm:gensall}
The left ideal $\genksn \twisted \lie$ is equal to $\liealg$.
\end{theorem}

By extending Definition \ref{def:lieidempotent} in the obvious way, 
we define what it means for an element of $\genksn$ to be a Lie idempotent,
and we have the following immediate
consequence of Theorems \ref{thm:idem} and \ref{thm:gensall}:

\begin{corollary}
$\lie$ is a Lie idempotent.
\end{corollary}

One might wonder why $\lie$ would have such desirable properties.  In Section 5, 
we give one possible explanation.  Removing the condition that $\ali_1\cdots \ali_n=1$, 
we consider an expression $\infproda$, which one can think of as a generating function for $\lie$, 
defined by
\[
\infproda = \sum_{n \geq 0} \frac{\e{n}{\ali_1, \ldots, \ali_n}}{(1-\ali_1\cdots \ali_n)}.
\]
As our main result of Section 5, we show that $\infproda$ can be expressed as a 
very simple infinite product.  This result generalizes \cite[Proposition 5.10]{GKLLRT95},
which corresponds to the specialization $q_i \to q$ for $i=1, \ldots, n$.  

\begin{remark}
According to the referee, our work possibly has a generalization in the spirit of the
papers of Lascoux, Leclerc and Thibon \cite{LLT95}, and of Hivert \cite{Hiv00}.
In \cite{LLT95}, a multi-parameter construction is devised, not for the Klyachko
idempotent, but for rectangular-shaped $q$-Kostka numbers.  
On the other hand, 
$q$-Kostka numbers are defined in terms of Hall-Littlewood polynomials, and \cite{Hiv00}
shows a direct connection between column-shaped
Hall-Littlewood polynomials and the Klyachko idempotent. 
\end{remark} 

The organization of the remainder of the paper is simple: in Sections 2, 3 and 4,
we prove Theorems \ref{thm:lie}, \ref{thm:idem} and \ref{thm:gensall} respectively.  
The infinite product expansion is the subject of Section 5. 

Before beginning the proofs, we need to introduce some terminology related to 
permutations.
If $w$ is a word of length $n$, we will write $w(i)$ to denote the $i$th
letter
of $w$.  
If the letters of $w$ are distinct, we define the \emph{descent set} $\des{w}$
of $w$ by 
$\des{w} = \{ i\ |\ 1 \leq i \leq n-1, w(i) > w(i+1)\}$.  
The \emph{major index} $\maj{w}$ of $w$ is then the sum of 
the elements of $\des{w}$.  
We will denote the cardinality
of $\des{w}$ by $d(w)$. 
Finally, we will use $\bar{D}(w)$ to denote the of \emph{circular descent set} of $w$, so 
that 
\[
\bar{D}(w) = \left\{ \begin{array}{ll}
D(w) & \mbox{if $w(n) < w(1)$} \\
D(w) + \{n\}  & \mbox{if $w(n) > w(1)$}
\end{array} \right. .
\] 
Then $\bar{d}(w)$ is simply the cardinality of $\bar{D}(w)$.  


\section{$\lie$ is a Lie element}

Our goal for this section is to prove Theorem \ref{thm:lie}.
We begin by stating a well-known characterization of Lie elements.
We refer the reader to \cite{Gar90}, \cite[p. 87]{Lot83} or  \cite[\S\S 1.3-1.4]{Reu93} for further details.

Our definitions in this paragraph will hold for $\mathbb{K}$ being any field of characteristic $0$,
although we will only need them for the case $\mathbb{K} = \genk$.  
First, define an inner (scalar) product $\scalar{\ }{\ }$ on 
$\mathbb{K}\langle X\rangle$ by $\scalar{u}{v} = \delta_{u,v}$ for any 
words $u$ and $v$, extended to $\mathbb{K}\langle X\rangle$ by linearity.
For the remainder of this section, it suffices to restrict to the case when $X = \{1, 2, \ldots, n\}$.   
Suppose
$u=u(1) u(2) \ldots u(r)$ and $v = v(1) v(2) \ldots v(s)$ and, since it will be sufficiently
general for our needs, assume
$u(i) \neq v(j)$ for all $i, j$.  A word $w$ is said to be a \emph{shuffle} of $u$ and $v$
if $w$ has length $r+s$ and if $u$ and $v$ are both subsequences of $w$.  The 
\emph{shuffle product} $u \shuffle v$ of $u$ and $v$ is an element of $\mathbb{K}\langle X \rangle$ 
and is defined to 
be the sum of all the shuffles of $u$ and $v$.  We will write $w \in u \shuffle v$
if $w$ is a shuffle of $u$ and $v$. 
The characterization of Lie elements that we will use is the following:
an element $\p$ of $\mathbb{K}\sn$ is a Lie element if and only if $\p$ is orthogonal to $u \shuffle v$ for
all non-empty words $u$ and $v$.  
Therefore, we wish to show that
\begin{equation}\label{equ:scalarzero}
\scalar{\lie}{u \shuffle v} = 0
\end{equation}
for all $u=u(1) u(2) \cdots u(r)$ and $v = v(1) v(2) \cdots v(s)$ with $r, s \geq 1$.  

Because \eqref{equ:scalarzero} holds trivially otherwise, let us assume that $r+s=n$ and that
$u(1), u(2), \ldots, u(r), v(1), v(2), \ldots, v(s)$ are all distinct.
Therefore, the partially ordered set (poset) $\puv$ whose Hasse
diagram is shown in Figure \ref{fig:puv} is a poset with elements $\{1,2,\ldots, n\}$.
Namely, $\puv$ is the disjoint union of the chains $u(1) < u(2) < \cdots < u(r)$ and
$v(1) < v(2) < \cdots < v(s)$.  
\setlength{\unitlength}{2pt}
\begin{figure} 
\begin{center}
\begin{picture}(40,45)
\thicklines
\put(0,0){\begin{picture}(20,45)(0,4)
\multiput(10,5)(0,10){2}{\circle*{2}}
\put(10,5){\line(0,1){15}}
\multiput(10,21)(0,2){10}{\line(0,1){1}}
\put(10,40){\line(0,1){5}}
\put(10,45){\circle*{2}}
\put(-1,4){$u(1)$}
\put(-1,14){$u(2)$}
\put(-1,44){$u(r)$}
\end{picture}}
\put(20,0){\begin{picture}(20,45)(0,4)
\multiput(10,10)(0,10){2}{\circle*{2}}
\put(10,10){\line(0,1){15}}
\multiput(10,26)(0,2){5}{\line(0,1){1}}
\put(10,35){\line(0,1){5}}
\put(10,40){\circle*{2}}
\put(-1,9){$v(1)$}
\put(-1,19){$v(2)$}
\put(-1,39){$v(s)$}
\end{picture}}
\end{picture}
\end{center}
\caption{$\puv$}
\label{fig:puv}
\end{figure}
\setlength{\unitlength}{1pt}
Recall that a \emph{linear extension} $\sigma$ of a poset $P$ of size $n$ is a bijection 
$\sigma:P \to \{1,2, \ldots, n\}$ such that if $y \leq z$ in $P$, then $\sigma(y) \leq \sigma(z)$.
We will represent the linear extension $\sigma$ as the word
$\sigma^{-1}(1), \sigma^{-1}(2), \ldots, \sigma^{-1}(n)$, and we will write
$\linext{P}$ to denote the set of linear extensions of $P$.  
We introduce linear extensions and the poset $\puv$ for the following reason:
the set of shuffles of $u$ and $v$ is exactly the set of linear extensions of $\puv$.
Therefore, $\scalar{\lie}{u \shuffle v}$ can be expressed as
\begin{equation}\label{equ:resemblance}
\scalar{\lie}{u \shuffle v}
= \sum_{\sigma \in \linext{\puv}} \frac{\prod_{j \in \des{\sigma}} \ali_{\sigma(1)} \ali_{\sigma(2)} \ldots \ali_{\sigma(j)}}{\prod_{i=1}^{n-1} (1 - \ali_{\sigma(1)} \ali_{\sigma(2)} \ldots \ali_{\sigma(i)})} .
\end{equation}
The reader who is acquainted with Richard Stanley's theory of $P$-partitions may find \eqref{equ:resemblance} strikingly familiar.  We now introduce the parts of this theory that
will be necessary to complete our proof.  While $P$-partitions are the topic of 
\cite[\S 4.5]{ec1}, we will need to work in the slightly more general setting found
in \cite{StaThesis}.

For our purposes, it is most convenient to say that a \emph{labelling} $\omega$ of a poset $P$
is an injection $\omega:P \hookrightarrow \{1, 2, \ldots\}$.
\begin{definition}
Let $P$ be a finite partially ordered set with a labelling $\omega$.  A 
\emph{$(P,\omega)$-partition} 
is a map $f:P \to \{0,1,2,\ldots\}$ with the following properties:
\begin{itemize}
\item[(i)]  $f$ is \emph{order-reversing}: if $y \leq z$ in $P$ then $f(y) \geq f(z)$,
\item[(ii)] if $y < z$ in $P$ and $\omega(y) > \omega(z)$, then $f(y) > f(z)$.
\end{itemize}
\end{definition}

In short, $(P,\omega)$-partitions are order-reversing maps with certain strictness conditions
determined by $\omega$.  We will denote the set of $(P,\omega)$-partitions by $\apw$.

\begin{note}
If $P$ is a poset with elements contained in the set $\{1,2,\ldots,n\}$, 
then in the labelled poset $(P,\omega)$, 
each vertex $i$ will have a label $\omega(i)$ associated to it and, in general,
we certainly
need not have $\omega(i) = i$.  However, in our case, we will always take
$\omega(i)=i$, since this is sufficient to yield the desired outcome.  
\end{note}

Define the generating function $\fpwx$ in the variables $\x = (x_1, x_2, \ldots, x_n)$ by
\[
\fpwx = \sum_{f \in \apw} \left( \prod_{p \in P} (x_{\omega(p)})^{f(p)} \right) .
\]
For any $n$-element poset $P$, 
by \cite[Prop. 7.1]{StaThesis} in the case when $\omega=\id$, the identity map, we have
\[
\fpwx = \sum_{\sigma \in \linext{P}} \frac{\prod_{j \in \des{\sigma}} x_{\sigma(1)} x_{\sigma(2)} \ldots x_{\sigma(j)}}{\prod_{i=1}^n (1 - x_{\sigma(1)} x_{\sigma(2)} \ldots x_{\sigma(i)})} .
\]
Comparing this with \eqref{equ:resemblance}, we deduce that
\[
\scalar{\lie}{u \shuffle v} = (1-\ali_1 \ali_2 \cdots \ali_n) F(\puv, \id; \al) .
\]
The structure of $\puv$ is simple enough that we can actually get a nice expression
for $F(\puv, \id; \al)$.  
Indeed, when $P$ is simply a total order with elements labelled 
$u(1), u(2), \ldots, u(r)$ from
bottom to top, we see that
\begin{equation}\label{equ:chain}
\fpwx = \frac{\prod_{j \in \des{u}} x_{u(1)} x_{u(2)} \ldots x_{u(j)}}
{\prod_{i=1}^r (1 - x_{u(1)} x_{u(2)} \ldots x_{u(i)})} ,
\end{equation}
where $u$ is the word $u(1) u(2) \ldots u(r)$.  
The terms in the denominator ensure that the $(P,\omega)$-partitions are order-reversing, 
while the terms in the numerator take care of the strictness conditions.  
Furthermore,
if $P$ is a disjoint union $P = P_1 + P_2$, then let $\omega_i$ denote the labelling
$\omega$ restricted to the elements of $P_i$, for $i=1,2$.  We see that
\begin{equation}\label{equ:disjointunion}
\fpwx = F(P_1, \omega_1; \x) F(P_2, \omega_2; \x).
\end{equation}
Combining \eqref{equ:chain} and \eqref{equ:disjointunion}, 
we deduce that
\[
F(\puv, \id; \x) = \frac{\left( \prod_{j \in \des{u}} x_{u(1)} x_{u(2)} \cdots x_{u(j)} \right) 
\left( \prod_{\ell \in \des{v}} x_{v(1)} x_{v(2)} \cdots x_{v(\ell)} \right)}
{\prod_{i=1}^r \left(1 - x_{u(1)} x_{u(2)} \cdots x_{u(i)} \right)
\prod_{k=1}^s \left(1 - x_{v(1)} x_{v(2)} \cdots x_{v(k)} \right) }.
\]
We finally conclude that 
\begin{eqnarray*}
\scalar{\lie}{u \shuffle v} & = & (1-\ali_1 \cdots \ali_n) F(\puv, \id; \al) \\
& = & (1-\ali_1 \cdots \ali_n)
\frac{\left( \prod_{j \in \des{u}} \ali_{u(1)} \cdots \ali_{u(j)} \right) 
\left( \prod_{\ell \in \des{v}} \ali_{v(1)} \cdots \ali_{v(\ell)} \right)}
{\left( \prod_{i=1}^r 1 - \ali_{u(1)} \cdots \ali_{u(i)} \right)
\left( \prod_{k=1}^s 1 - \ali_{v(1)} \cdots \ali_{v(k)} \right)}  \\
& = & 0 ,
\end{eqnarray*}
because of the conditions on $\al$ and because $r, s < n$. 
This yields Theorem \ref{thm:lie}.


\section{$\lie$ is idempotent}

One way to show that the Klyachko idempotent is idempotent 
is to define an element $\klypartner$ of $\ksn$ such that
\[
\klypartner \kly = \kly \mbox{\ \ and \ \ } \kly \klypartner = \klypartner.
\]
(See \cite{Kly74}, \cite[Lemma 8.19]{Reu93}.)
Then it follows that 
\[
\kly^2 = \kly \klypartner \kly = \klypartner \kly = \kly, 
\]
as required.  
Throughout, let $\cycle$ denote the $n$-cycle $(1,2,\ldots, n) \in \sn$.  
Let $\ru$ denote the primitive $n$th root of unity from the definition 
of $\kly$.
Then a suitable element $\klypartner$ is 
given by 
\[
\klypartner = \frac{1}{n} \sum_{i=0}^{n-1} \frac{\cycle^i}{\ru^i}.
\]

We wish to apply the same principle to show that $\lie$ is idempotent, thus proving Theorem \ref{thm:idem}.
We define an element $\prt$ of $\genksn$ by
\[
\prt = \sum_{i=0}^{n-1} \gmaja{\cycle^i} \cycle^i.
\]
The reader is encouraged to check that if
$\ali_1, \ldots, \ali_n$ are all mapped to $\ru$, 
then $\prt$ maps to $\klypartner$.
Our goal, therefore, for the remainder of this section is to show that:
\begin{eqnarray}
\prt \twisted \lie & = & \lie, \mbox{\ \ and}   \label{equ:idem_a} \\
\lie \twisted \prt & = & \prt.  \label{equ:idem_b}
\end{eqnarray} 

Because we are taking twisted products, we will need to know how, for example,
$\apply{\cycle}{\gmaja{\sigma}}$ compares to $\gmaja{\sigma}$.
For notational convenience, for any $\sigma \in \sn$, let us write 
$\gmaja{\sigma} = \frac{\N{\sigma}}{\D{\sigma}}$, with
\begin{eqnarray*}
\N{\sigma} & = & \prod_{j \in \des{\sigma}} \ali_{\sigma(1)} \ali_{\sigma(2)} \ldots \ali_{\sigma(j)},\\
\D{\sigma} & = & \prod_{i=1}^{n-1} (1 - \ali_{\sigma(1)} \ali_{\sigma(2)} \ldots \ali_{\sigma(i)}).
\end{eqnarray*}
(``N'' stands for numerator, and ``D'' for denominator.)
The following result extends \cite[Lemma 11]{PaRe02}.

\begin{lemma}\label{lem:applygamma}
For all $\sigma, \tau \in \sn$, we have the following identities among elements of $\genk$:
\begin{enumerate}
\item[(i)] \[  \apply{\cycle}{\N{\sigma}} = \ali_1 \N{\cycle \sigma}.\]
\item[(ii)] \[ \apply{\tau}{\D{\sigma}} = \D{\tau \sigma}. \]
\item[(iii)] \[ \apply{\cycle^i}{\gmaja{\sigma}}   =  \ali_1 \cdots \ali_i \cdot \gmaja{\cycle^i \sigma}.\]
\item[(iv)] \[ \N{\sigma \cycle^i}  =  (\ali_{\sigma(1)}\cdots \ali_{\sigma(i)})^{-\bar{d}(\sigma)} \N{\sigma}.\]
\end{enumerate}
\end{lemma}

\begin{proof}
(i) By definition, 
\begin{eqnarray*}
\apply{\cycle}{\N{\sigma}}  & = & 
\prod_{j \in \des{\sigma}} \ali_{\cycle( \sigma(1))} \ldots \ali_{\cycle(\sigma(j))}, \\
& = & \prod_{j \in \des{\sigma}} \ali_{(\cycle \sigma)(1)} \ldots \ali_{(\cycle\sigma)(j)},
\end{eqnarray*}
since $\cycle(\sigma(i)) = (\cycle \sigma)(i)$.

The argument that follows is best understood by first trying some simple examples.
If $\sigma(n) = n$, then $\des{\cycle \sigma} = \des{\sigma} + \{n-1\}$.
Therefore, 
\begin{eqnarray*}
\N{\cycle \sigma} & = & (\ali_{\cycle\sigma(1)} \ali_{\cycle\sigma(2)} \cdots \ali_{\cycle\sigma(n-1)}) \apply{\cycle}{\N{\sigma}} \\
& = & (\ali_2 \ali_3 \cdots \ali_n) \apply{\cycle}{\N{\sigma}}\\
& = & \frac{\apply{\cycle}{\N{\sigma}}}{\ali_1}.
\end{eqnarray*}
If $\sigma^{-1}(n) = i < n$, then $\des{\cycle \sigma} = \des{\sigma} + \{i-1\} - \{i\}$,
where we set $\{0\} = \emptyset$.  Therefore, 
\[
\N{\cycle \sigma} = \frac{\ali_{\cycle\sigma(1)} \cdots \ali_{\cycle\sigma(i-1)}}
{\ali_{\cycle\sigma(1)} \cdots \ali_{\cycle\sigma(i)}}
\apply{\cycle}{\N{\sigma}} = \frac{\apply{\cycle}{\N{\sigma}}}{\ali_1}.
\]

(ii) This follows directly from the fact that $\tau(\sigma(i)) = (\tau \sigma)(i)$.   

(iii) By (i) and (ii), this is clearly true when $i=1$.  Working by induction, 
\begin{eqnarray*}
\apply{\cycle^i}{\gmaja{\sigma}} & = & \apply{\cycle}{\apply{\cycle^{i-1}}{\gmaja{\sigma}}} \\
& = & \apply{\cycle}{\ali_1 \cdots \ali_{i-1} \cdot \gmaja{\cycle^{i-1}\sigma}} \\
& = & \ali_2 \cdots \ali_i \cdot \apply{\cycle}{\gmaja{\cycle^{i-1}\sigma}} \\
& = & \ali_1 \ali_2 \cdots \ali_i \cdot \gmaja{\cycle^i \sigma}.
\end{eqnarray*}

(iv) We first show that 
\begin{equation}\label{equ:nsigma}
\N{\sigma \cycle} = \frac{\N{\sigma}}{(\ali_{\sigma(1)})^{\bar{d}(\sigma)}}.
\end{equation}
Indeed, suppose that $\sigma(1) > \sigma(n)$.  Then
\[
\N{\sigma \cycle} = \frac{\N{\sigma}}{(\ali_{\sigma(1)})^{d(\sigma)}}
\]
holds directly, and $d(\sigma) = \bar{d}(\sigma)$.
If $\sigma(1) < \sigma(n)$, then
\begin{eqnarray*}
\N{\sigma \cycle} & = & \frac{\ali_{\sigma(2)} \cdots \ali_{\sigma(n)}\N{\sigma}}
{(\ali_{\sigma(1)})^{d(\sigma)}}\\
& = & \frac{\N{\sigma}}{(\ali_{\sigma(1)})^{d(\sigma)+1}},
\end{eqnarray*}
and $d(\sigma)+ 1 = \bar{d}(\sigma)$.

Proceeding by induction, 
\begin{eqnarray*}
\N{\sigma \cycle^i}  & =  & \N{(\sigma \cycle^{i-1}) \cycle} \\
& = & \frac{\N{\sigma \cycle^{i-1}}}{(\ali_{\sigma \cycle^{i-1}(1)})^{\bar{d}(\sigma \cycle^{i-1})}} \\
& = & \frac{\N{\sigma \cycle^{i-1}}}{(\ali_{\sigma(i)})^{\bar{d}(\sigma)}} \\
& = & (\ali_{\sigma(1)}\cdots \ali_{\sigma(i)})^{-\bar{d}(\sigma)} \N{\sigma}.
\end{eqnarray*}

\end{proof}

Before proving \eqref{equ:idem_a} and \eqref{equ:idem_b}, we state one further
necessary result, which is essentially taken word-for-word from \cite{PaRe02}.
As usual, $\delta_{i,j}$ denotes the Kronecker delta, defined to be 1
if $i=j$, and 0 otherwise.

\begin{proposition}\label{pro:pare} \cite[Corollary 10]{PaRe02}
Let $\alpha_1, \ldots, \alpha_n$ be elements of a field such that 
the product $\alpha_1 \alpha_2 \cdots \alpha_n$ 
is equal to 1 and each subproduct is different from 1.  For
$k = 0, 1, \ldots, n-1,$
\[
\sum_{i = 0}^{n-1} \frac{(\alpha_{i+1} \ldots \alpha_n)^k}
{(1 - \alpha_{i+1})(1 - \alpha_{i+1}\alpha_{i+2}) \cdots (1 - \alpha_{i+1}\alpha_{i+2} \cdots \alpha_{i+n-1})} = \delta_{0,k},
\]
where all subscripts on $\alpha$ are taken modulo $n$. 
\end{proposition}

Here and henceforth, unless otherwise stated, all subscripts on $\al$ are taken
modulo $n$. 

\begin{proof}[Proof of \eqref{equ:idem_a}]
We wish to show that
\begin{equation}\label{equ:idem_a_full}
\left( \sum_{i=0}^{n-1} \gmaja{\cycle^i} \cycle^i \right) \twisted 
\left( \sum_{\sigma \in \sn} \gmaja{\sigma} \sigma \right) = 
 \sum_{\tau \in \sn} \gmaja{\tau} \tau.
\end{equation}
The left-hand side can be rewritten as 
\begin{eqnarray*}
& & \sum_{i, \sigma} \gmaja{\cycle^i} \cdot \apply{\cycle^i}{\gmaja{\sigma}} \cycle^i \sigma \\
& = & \sum_{i, \sigma} \gmaja{\cycle^i} \cdot \ali_1 \cdots \ali_i \cdot 
\gmaja{\cycle^i \sigma} \cycle^i \sigma \\
\end{eqnarray*}
by Lemma \ref{lem:applygamma}(iii).  Therefore, showing \eqref{equ:idem_a_full} is equivalent
to showing that
\[
\sum_{\genfrac{}{}{0pt}{}{i, \sigma}{\cycle^i \sigma = \tau}}
\gmaja{\cycle^i} \cdot \ali_1 \cdots \ali_i \cdot \gmaja{\cycle^i \sigma} 
= \gmaja{\tau}
\]
for all $\tau \in \sn$.  
Here, the left-hand side simplifies to
\begin{eqnarray}
& & \gmaja{\tau} \sum_{\genfrac{}{}{0pt}{}{i, \sigma}{\cycle^i \sigma = \tau}}
\gmaja{\cycle^i} \cdot \ali_1 \cdots \ali_i \nonumber \\
& = & \gmaja{\tau} \sum_{i=0}^{n-1} \gmaja{\cycle^i} \cdot \ali_1 \cdots \ali_i \nonumber \\
& = & \gmaja{\tau} \sum_{i=0}^{n-1} \frac{1}{(1-\ali_{i+1})(1-\ali_{i+1}\ali_{i+2})
\cdots (1 - \ali_{i+1}\ali_{i+2} \cdots \ali_{i+n-1})} \label{equ:lonelabel}
\end{eqnarray}
by Lemma \ref{lem:applygamma}(iii) with $\sigma = \id$ or, alternatively, by direct calculation
of  $\gmaja{\cycle^i}$.  Applying Proposition \ref{pro:pare}, we see that the expression
of \eqref{equ:lonelabel} equals  $\gmaja{\tau}$, as required.
\end{proof}

\begin{proof}[Proof of \eqref{equ:idem_b}]
We wish to show that
\begin{equation}\label{equ:idem_b_full}
\left( \sum_{\sigma \in \sn} \gmaja{\sigma} \sigma \right)  \twisted 
\left( \sum_{i=0}^{n-1} \gmaja{\cycle^i} \cycle^i \right) 
=  \sum_{j=0}^{n-1} \gmaja{\cycle^j} \cycle^j .
\end{equation}
Since the left-hand side can be rewritten as
\[
\sum_{\sigma, i}
\gmaja{\sigma} \cdot \apply{\sigma}{\gmaja{\cycle^i}} \sigma \cycle^i ,
\]
showing \eqref{equ:idem_b_full} is equivalent to showing that
\begin{equation}\label{equ:twopossibilities}
\sum_{\genfrac{}{}{0pt}{}{\sigma, i}{\sigma \cycle^i = \tau}}
\gmaja{\sigma} \cdot \apply{\sigma}{\gmaja{\cycle^i}}
= \left\{ \begin{array}{ll} \gmaja{\cycle^j}  & \mbox{if $\tau = \cycle^j$ for some $j$,} \\
0 & \mbox{otherwise.} \end{array} \right.
\end{equation}
The left-hand side of \eqref{equ:twopossibilities} can be rewritten as
\begin{eqnarray*}
& & \sum_{i=0}^{n-1} \gmaja{\tau \cycle^{-i}} \cdot \apply{(\tau \cycle^{-i})}{\gmaja{\cycle^i}} \\
& = & \sum_{i=0}^{n-1} \gmaja{\tau \cycle^i} \cdot \apply{(\tau \cycle^i)}{\gmaja{\cycle^{n-i}}} .
\end{eqnarray*}
By Lemma \ref{lem:applygamma}(iv), 
\[ 
\N{\tau \cycle^i}  =  (\ali_{\tau(1)}\cdots \ali_{\tau(i)})^{-\bar{d}(\tau)} \N{\tau}.
\]
Also, we see directly that
\[
\D{\tau \cycle^i} = (1-\ali_{\tau(i+1)})(1-\ali_{\tau(i+1)}\ali_{\tau(i+2)}) \cdots 
(1-\ali_{\tau(i+1)} \ali_{\tau(i+2)} \cdots \ali_{\tau(i+n-1)}).
\]
As for the term $\apply{(\tau \cycle^i)}{\gmaja{\cycle^{n-i}}}$, we have that
\begin{eqnarray*}
\apply{(\tau \cycle^i)}{\N{\cycle^{n-i}}} & = & \apply{(\tau \cycle^i)}{\ali_{n-i+1} \cdots \ali_n} \\
& = & \apply{\tau}{\ali_1 \cdots \ali_i} \\
& = & \ali_{\tau(1)} \cdots \ali_{\tau(i)},
\end{eqnarray*}
while, by Lemma \ref{lem:applygamma}(ii), 
$\apply{(\tau \cycle^i)}{\D{\cycle^{-i}}} = \D{\tau}$.  

Putting this all together, we get that
\begin{eqnarray*}
& & \sum_{\genfrac{}{}{0pt}{}{\sigma, i}{\sigma \cycle^i = \tau}} 
\gmaja{\sigma} \cdot \apply{\sigma}{\gmaja{\cycle^i}} \\
& = &
\frac{\N{\tau}}{\D{\tau}} \sum_{i=0}^{n-1} \frac{(\ali_{\tau(1)}\cdots \ali_{\tau(i)})^{1-\bar{d}(\tau)}}
{(1-\ali_{\tau(i+1)})
\cdots (1-\ali_{\tau(i+1)} \ali_{\tau(i+2)} \cdots \ali_{\tau(i+n-1)})} \\
& = &
\frac{\N{\tau}}{\D{\tau}} \sum_{i=0}^{n-1} \frac{(\ali_{\tau(i+1)}\cdots \ali_{\tau(n)})^{\bar{d}(\tau)-1}}
{(1-\ali_{\tau(i+1)})
\cdots (1-\ali_{\tau(i+1)} \ali_{\tau(i+2)} \cdots \ali_{\tau(i+n-1)})} \ .
\end{eqnarray*}
Since $\bar{d}(\tau) = 1$ if $\tau = \cycle^j$ for some $j$, and 
$2 \leq \bar{d}(\tau) \leq n-1$ otherwise, applying Proposition \ref{pro:pare} gives
exactly the desired equality \eqref{equ:twopossibilities}.  
\end{proof}


\section{$\lie$ generates the multilinear part of the free Lie algebra}

As before, let $\liealgn$ (resp. $\liealg$)
denote the set of Lie elements in $\ksn$ (resp. $\genksn$).
Our goal for this section is to show that
\[
\genksn \twisted \lie = \liealg.
\]
By Theorem \ref{thm:lie}, we know that $\lie \in \liealg$, while \eqref{equ:liealgleftideal}
states that $\liealg$ is a left ideal of $\genksn$.  Hence, 
we have that $\genksn \twisted \lie \subseteq \liealg$.
Also, it is well-known that $\liealgn$ has dimension $(n-1)!$, from which it follows that
$\liealg$ has dimension $(n-1)!$.
Therefore, it suffices to show that $\genksn \twisted \lie$ has dimension at least $(n-1)!$.  

As we saw in Section 3, we have the following two identities:
\begin{eqnarray*}
\prt \twisted \lie & = & \lie,  \\
\lie \twisted \prt & = & \prt.  
\end{eqnarray*} 
We claim that these identities can be used to give a bijection 
$\phi$ from $\genksn \twisted \lie$ to $\genksn \twisted \prt$.
Indeed, for $x \in \genksn \twisted \lie$, let 
\[
\phi(x) = x \twisted \prt.
\]
Define a map $\psi : \genksn \twisted \prt \to \genksn \twisted \lie$ by 
\[
\psi(y) = y \twisted \lie
\]
for any $y \in \genksn \twisted \prt$.  
For $x \in \genksn \twisted \lie$, we know that $x = x' \twisted \lie$ for some $x' \in \genksn$.  
Therefore, 
\begin{eqnarray*}
\psi(\phi(x)) & = & \psi(x' \twisted \lie \twisted \prt) \\
& = & \psi(x' \twisted \prt) \\
& = & (x' \twisted \prt) \twisted \lie \\
& = & x' \twisted \lie \\
& = & x.
\end{eqnarray*}
Similarly, $\phi(\psi(y)) = y$ for any $y \in \genksn \twisted \prt$, and so $\phi$ is a 
bijection.
We conclude that to prove Theorem \ref{thm:gensall}, it remains to show 
that $\genksn \twisted \prt$ has dimension at least $(n-1)!$.  

Consider the set $B = \{ \sigma \twisted \prt\ |\ \sigma\in \sn, \sigma(1)=1 \}$.
Clearly, $B \subseteq \genksn \twisted \prt$.  
We claim that $B$ forms a basis of $\genksn \twisted \prt$.  
Suppose $\sigma$ and $\tau$ are distinct permutations with $\sigma(1) = \tau(1) = 1$.
Since
\[
\sigma \twisted \prt = \sum_{i=0}^{n-1} \apply{\sigma}{\gmaja{\cycle^i}} \sigma \cycle^i  ,
\]
we see that $\sigma \twisted \prt$ is a linear combination of permutations of the form $\sigma \cycle^i$.
Similarly, $\tau \twisted \prt$ is a linear combination of permutations of the form $\tau \cycle^j$.
But since $\sigma$ and $\tau$ are distinct and $\sigma(1) = \tau(1) = 1$, there do
not exist $k, \ell$ such that $\sigma \cycle^k = \tau \cycle^\ell$.  It follows that 
$\sigma \twisted \prt \neq \tau \twisted \prt$ and, furthermore, that the elements of $B$ are
linearly independent.  We conclude that $B$ consists of $(n-1)!$ linearly
independent elements, and so $\genksn \twisted \prt$ has dimension at least
$(n-1)!$, thus proving Theorem \ref{thm:gensall}.

We will conclude by showing independently that $B$ spans $\genksn \twisted \prt$, thereby
reproving that the dimension of $\genksn \twisted \prt$, 
and hence of $\genksn \twisted \lie$, is $(n-1)!$.

Every permutation $\tau \in \sn$ is of the form $\sigma \cycle^j$ for some $\sigma$ with
$\sigma(1) = 1$.  
We claim that $\tau \twisted \prt$ is then simply a scalar multiple of $\sigma \twisted \prt \in B$.
Indeed, using Lemma \ref{lem:applygamma}(iii), we have
\begin{eqnarray*}
\tau \twisted \prt & = & (\sigma \cycle^j) \twisted \sum_{i=0}^{n-1} \gmaja{\cycle^i} \cycle^i \\
& = & \sum_{i=0}^{n-1} \apply{(\sigma \cycle^j)}{\gmaja{\cycle^i}} \sigma \cycle^{j+i} \\
& = &  \sum_{i=0}^{n-1} \ali_{\sigma(1)} \cdots \ali_{\sigma(j)} \cdot
\apply{\sigma}{\gmaja{\cycle^{j+i}}} \sigma \cycle^{j+i} \\
& = & \ali_{\sigma(1)} \cdots \ali_{\sigma(j)}  \sum_{i=0}^{n-1} 
\apply{\sigma}{\gmaja{\cycle^i }} \sigma \cycle^i \\
& = & \ali_{\sigma(1)} \cdots \ali_{\sigma(j)} \cdot \sigma \twisted \prt.
\end{eqnarray*}
It follows that every element of $\genksn \twisted \prt$ can be written as a linear combination
of elements of $B$, as required.


\section{Infinite product expansion}

Let $\x = \xj_1, \xj_2, \ldots$ be an infinite sequence of variables in a field $\K$ of characteristic 0.
In effect, the $\xj_i$'s will play the role formerly played by the $\ali_i$'s.  We switch variables
to emphasize two key differences with the material in this section.  The first is that $n$ will no
longer be fixed and the second is that there will be no restriction on the $\xj_i$'s analogous to
the restriction $\ali_1 \cdots \ali_n = 1$ on the $\ali_i$'s.  
We let $\kxn$ denote the algebra of 
formal power series in $\xj_1, \ldots, \xj_n$ over the field $\K$.

The central object of study in this section will be $\infprod$, which is an element of
$\directsum$, and is defined as follows:

\begin{eqnarray}\label{equ:thetadef}
\infprod & = & \sum_{n \geq 0} \sum_{\sigma \in \sn} \frac{
\prod_{j \in \des{\sigma}} \xj_{\sigma(1)} \xj_{\sigma(2)} \ldots \xj_{\sigma(j)}}
{\prod_{i=1}^n (1 - \xj_{\sigma(1)} \xj_{\sigma(2)} \ldots \xj_{\sigma(i)})} \sigma \\
& = & \sum_{n\geq 0} \frac{\e{n}{\xj_1, \ldots, \xj_n}}{(1-\xj_1 \cdots \xj_n)} \nonumber.
\end{eqnarray}

We can think of $\infprod$ as a generating function for $\e{n}{\xj_1, \ldots, \xj_n}$ and, in particular, 
we see that knowing $\infprod$ allows us to extract $\lie$, for any $n$.  

Our goal for this section is to show that $\infprod$ can be expressed as a simple
infinite product.  To do this, we first need to put an appropriate 
ring structure on 
$\directsum$.  As we shall see, once we have done this correctly, we will
have completed much of the work necessary to fulfill our goal.

First, we will recall from \cite{MaRe95} the associative product $\mare$ defined on 
$\oplus_{n\geq 0} \ksn$.  We let $\mathbb{P}$ denote the positive integers, and we will be 
working with words on the alphabet $\mathbb{P}$.  If $u$ and $v$ are two such words, 
then we define the product $u \cdot v$ to be simply their concatenation.  
If $w = w(1) \ldots w(n)$ and the letters $w(1), \ldots, w(n)$ of $w$ are 
distinct, then we define the \emph{standardization} of $w$ to be the unique
permutation $\st(w)$ of $\sn$ satisfying 
\[
\st(w)(i) \leq \st(w)(j) \   \Leftrightarrow\ w(i) \leq w(j)
\]
for all $1 \leq i,j \leq n$.  For example, $\st(5716) = 2413$.
If $\sigma \in S_s$ and $\tau \in S_t$ then we define $\mare$
by 
\[
\sigma \mare \tau = \sum_{u,v} u \cdot v ,
\]
where the sum is over all $u, v$ such that $u \cdot v \in S_{s+t}$ with $\st(u) = \sigma$ and
$\st(v) = \tau$.  For example, $132 \mare 1 = 1324 + 1423 + 1432 + 2431$.  
Extending by linearity gives a product on $\oplus_{n\geq 0} \ksn$.
Since $\sigma \mare \tau$ is a multiplicity-free sum, we can write $w \in \sigma \mare \tau$ to
mean that $w$ appears as a term in $\sigma \mare \tau$.  

Before extending $\mare$ to $\directsum$, we will introduce some convenient notation.
If $f(\xj_1, \ldots, \xj_n) \in \kxn$ and if $w$ is a word on $\mathbb{P}$
of length $n$ with distinct letters $a_1 < a_2 < \cdots < a_n$,
then define $f(w)$ to be $f(\xj_{a_1}, \ldots, \xj_{a_n})$.  For example, 
if $f(\xj_1, \xj_2, \xj_3) \in \K[[\xj_1, \xj_2, \xj_3]]$, then  $f(382) = f(\xj_2, \xj_3, \xj_8)$.  
For $\sigma \in S_s$ and $\tau \in S_t$, we can then define
\[
(f(\xj_1, \ldots, \xj_s) \sigma) \mare (g(\xj_1, \ldots, \xj_t) \tau)
= \sum_{w \in \sigma \mare \tau} f(w(1) \ldots w(s)) \cdot g(w(s+1) \ldots w(s+t)) \, w, 
\]
and extend to $\directsum$ by linearity.  
For example, 
\begin{eqnarray*}
(\xj_1^2 \xj_2^2 \xj_3 132) \mare (\xj_1^3 1) & = & \xj_1^2 \xj_2^2 \xj_3 \xj_4^3 \, 1324 +
\xj_1^2 \xj_2^2 \xj_4 \xj_3^3 \, 1423 \\
& &  + \xj_1^2 \xj_3^2 \xj_4 \xj_2^3\, 1432 + \xj_2^2 \xj_3^2 \xj_4 \xj_1^3\, 2431.
\end{eqnarray*}
We see that the product $\mare$ on $\directsum$ is a natural extension of the 
$\oplus_{n\geq 0} \ksn$ version.  

\begin{remark}\label{rem:hopf}
The product $\mare$ on $\directsum$ makes it an associative algebra.  This follows from
\cite[\S 1.2]{BaHo05pr}: our algebra is the completion, in an appropriate sense, 
of the algebra $\oplus_{n\geq 0}\K[\xj_1, \ldots, \xj_n] \sn$;
the latter is isomorphic to the $\mathscr{F}(V)$ of Baumann and Hohlweg when
$V$ is taken to be the free $\K$-module with basis $\mathbb{N}$ (the non-negative integers).  
Indeed, a word 
$i_1 \ldots i_n$ on $\mathbb{N}$ corresponds to the monomial $\xj_1^{i_1} \cdots \xj_n^{i_n}$
in $\kxn$.  It also follows from \cite{BaHo05pr} that $\directsum$ can even be endowed with
a Hopf algebra structure.  
\end{remark}

We are almost ready to state the main result of this section.  We let $\emptyperm$ denote the
unique element of the symmetric group $S_0$.  When expressing an infinite product of 
elements of $\directsum$, we will use an arrow pointing left above
the product symbol to denote that the product should be expanded from right to left.  For example,
\[
\fancyprod{n\geq 0} F(n) = \cdots \mare F(2) \mare F(1) \mare F(0).
\]
The following result extends \cite[Proposition 5.10]{GKLLRT95}.

\begin{theorem} The generating function $\infprod$ of $\e{n}{\xj_1, \ldots, \xj_n}$ satisfies
\begin{equation}\label{equ:infprodexp}
\infprod = \fancyprod{n \geq 0} 
\left( \emptyperm + \xj_1^n\, 1 + \xj_1^n \xj_2^n\, 12 + \xj_1^n \xj_2^n \xj_3^n\, 123 + \cdots \right).
\end{equation}
\end{theorem}

The simplicity of \eqref{equ:infprodexp} perhaps helps to explain why $\lie$ would have many 
nice properties, such as those from the earlier sections.  We remark that each factor
in the product, and the product itself, is a group-like element for the Hopf algebra
structure of Remark \ref{rem:hopf}.  This is similar to 
the situation in \cite{GKLLRT95} and is easily verified.

\begin{proof}
We can expand the right-hand side of \eqref{equ:infprodexp} into nested sums as follows:
\begin{eqnarray*}
& & \fancyprod{n \geq 0}
\left( \emptyperm + \xj_1^n\, 1 + \xj_1^n \xj_2^n\, 12 + \xj_1^n \xj_2^n \xj_3^n\, 123 + \cdots \right) \\
& = & \sum_{k \geq 0} \ \sum_{n_1 > \cdots > n_k \geq 0} \ \sum_{a_1, \ldots, a_k \geq 1}
(\xj_1^{n_1} \cdots \xj_{a_1}^{n_1}\, 12\cdots a_1) \mare \cdots \mare 
(\xj_1^{n_k} \cdots \xj_{a_k}^{n_k}\, 12\cdots a_k) \\
& = & \sum_{n\geq 0} \ \sum_{\sigma \in \sn} \ \sum_{k \geq 0}\ \sum_{\sigma=u_1 \cdots u_k}\ 
\sum_{n_1 > \cdots > n_k \geq 0}  (\xj_{u_1})^{n_1}\cdots (\xj_{u_k})^{n_k}\, \sigma .
\end{eqnarray*}
where the fourth sum is over all non-empty increasing words $u_1, \ldots, u_k$ 
in $\mathbb{P}$ whose 
concatenation is the word $\sigma$, 
and where $\xj_u$ for the word
$u = u(1) \ldots u(s)$ denotes $\xj_{u(1)}\cdots \xj_{u(s)}$.  

We see that the coefficient $C(\sigma)$ of a fixed $\sigma \in \sn$ is then given by
\[
C(\sigma) = \sum_{k \geq 0}\ \sum_{\sigma=u_1\cdots u_k}\ 
\sum_{n_1 > \cdots > n_k \geq 0}  (\xj_{u_1})^{n_1}\cdots (\xj_{u_k})^{n_k} ,
\]
where the $u_i$'s are increasing non-empty words.
It will be helpful to rewrite $C(\sigma)$ as
\begin{eqnarray*}
C(\sigma) & = & \sum_{k \geq 0}\ \sum_{\sigma=u_1\cdots u_k}\ 
\sum_{n_1 > \cdots > n_k \geq 0}  (\xj_{u_1})^{n_1-n_2}
(\xj_{u_1 u_2})^{n_2 - n_3}
\cdots (\xj_{u_1 \ldots u_k})^{n_k} \\
& = & \sum_{k \geq 0}\ \sum_{\sigma=u_1\cdots u_k} \ 
\sum_{\genfrac{}{}{0pt}{}{m_1, \ldots, m_{k-1} > 0}{m_k \geq 0}}  (\xj_{u_1})^{m_1}
(\xj_{u_1 u_2})^{m_2}
\cdots (\xj_{u_1 \ldots u_k})^{m_k} .
\end{eqnarray*}
The condition that the $u_i$'s be increasing words implies that
$C(\sigma)$ can be restructured as
\[
C(\sigma) = \sum_{\genfrac{}{}{0pt}{}{p_1, p_2 \ldots, p_n \geq 0}
{i \in D(\sigma)\, \Rightarrow\, p_i > 0}}
 (\xj_{\sigma(1)})^{p_1}
(\xj_{\sigma(1)}\xj_{\sigma(2)})^{p_2}
\cdots
(\xj_{\sigma(1)}\cdots \xj_{\sigma(n)})^{p_n}.
\]
Observing that each monomial in the sum has 
$\prod_{j \in \des{\sigma}} \xj_{\sigma(1)} \xj_{\sigma(2)} \ldots \xj_{\sigma(j)}$ as a factor, we get
\begin{eqnarray}
C(\sigma) & = & \left( \prod_{j \in \des{\sigma}} \xj_{\sigma(1)} \xj_{\sigma(2)} \ldots \xj_{\sigma(j)} \right)
\times \\
& & 
\sum_{p_1, p_2 \ldots, p_n \geq 0} (\xj_{\sigma(1)})^{p_1} (\xj_{\sigma(1)}\xj_{\sigma(2)})^{p_2}
\cdots
(\xj_{\sigma(1)}\cdots \xj_{\sigma(n)})^{p_n}
 \nonumber \\
& = & \frac{\prod_{j \in \des{\sigma}} \xj_{\sigma(1)} \xj_{\sigma(2)} \ldots \xj_{\sigma(j)}}
{\prod_{i=1}^n (1 - \xj_{\sigma(1)} \xj_{\sigma(2)} \ldots \xj_{\sigma(i)})}, \label{equ:picompatible}
\end{eqnarray}
as required.
\end{proof}

We close with a remark that again relates our work to Stanley's $P$-partitions.  Comparing
\eqref{equ:picompatible} with \cite[Lemma 4.5.2(a)]{ec1}, we see that $C(\sigma)$ is exactly the generating function 
for the set of all $\sigma$-compatible permutations, as defined in \cite{ec1}.


\section*{Acknowledgements}
We thank Fran{\c{c}}ois Bergeron for suggesting that we look for an
infinite product expansion generalizing \cite[Proposition 5.10]{GKLLRT95},
and Pierre Baumann
for showing us 
the connection between $\directsum$ and the $\mathscr{F}(V)$ of \cite{BaHo05pr}, 
as discussed in Remark \ref{rem:hopf}.
Also, we are grateful for the referee's careful reading and suggestions.  


\bibliography{../../master}
\bibliographystyle{plain}

\end{document}